\documentclass[11pt,leqno]{article}

\usepackage{amsmath,amssymb,amsthm,mathrsfs,dsfont}

\usepackage[margin=3cm]{geometry} 

\usepackage{titlesec,hyperref}

\usepackage{color}

\usepackage{fancyhdr}
\usepackage{subfigure}
\usepackage[graphicx]{realboxes}
\usepackage{float}
\usepackage{enumerate}
\usepackage{enumitem}
\setlist[enumerate]{label=(\arabic*), leftmargin=*, itemsep=2pt}
\usepackage{dsfont} 
\usepackage{extarrows}

\titleformat{\subsection}{\it}{\thesubsection.\enspace}{1pt}{}

\newtheorem{theo}{Theorem}[section]
\newtheorem{lemm}[theo]{Lemma}

\newtheorem{prop}[theo]{Proposition}
\newtheorem{rema}[theo]{Remark}

\allowdisplaybreaks 

\newcommand{\R}{\mathbb R}
\newcommand{\A}{\mathcal A}\newcommand{\eps}{\varepsilon}
\newcommand{\dd}{\,\mathrm d}\newcommand{\norm}[1]{\left\|#1\right\|}
\DeclareMathOperator{\diver}{div}\DeclareMathOperator{\curl}{curl}

\numberwithin{equation}{section}

\begin{document}
	\title{Blow-up phenomena and ill-posedness for the 2D Euler-Poincar\'e equations
		\hspace{-4mm}
	}
	
	\author{
		Wei $\mbox{Luo}^1$ \footnote{E-mail:  luow23@mail.sysu.edu.cn}, \quad
		Zhaoyang $\mbox{Yin}^{1,2}$\footnote{E-mail: mcsyzy@mail.sysu.edu.cn} \quad and
		Pei $\mbox{Zheng}^{1}$\footnote{E-mail: zhengp25@mail2.sysu.edu.cn}\\
		$^1\mbox{Department}$ of Mathematics,
		Sun Yat-sen University, Guangzhou 510275, China\\
		$^2\mbox{School}$ of Science,\\ Shenzhen Campus of Sun Yat-sen University, Shenzhen 518107, China}
	
	\date{}
	\maketitle
	\hrule
	
	\begin{abstract}
		In this paper, we study finite-time blow-up and critical norm inflation for the Euler--Poincar\'e equations on $\mathbb R^2$. We consider a class of momentum fields admitting the div--curl decomposition with radial scalar potentials. The potentials satisfy a nonlocal transport system. By virtue of the conservation law and the finite propagation speed property, we obtain a growth estimate for the potential, which yields a sufficient condition for finite-time blow-up. Furthermore, combining a multiscale construction of the initial data with a blow-up criterion formulated in terms of the Besov norm of the deformation tensor, we establish norm inflation at the zero initial datum in $H^2(\mathbb R^2)$, thereby proving ill-posedness in this critical Sobolev space.
		
		\vspace*{5pt}
		
		\vspace*{5pt}
		
	\end{abstract}
	
	\vspace*{10pt}
	
	\tableofcontents

	\section{Introduction}
	In this paper, we consider the two-dimensional Euler--Poincar\'e equations
	\begin{equation}\label{EP}
		\left\{\begin{aligned}
			&\partial_t m+u\cdot\nabla m+(\nabla u)^T m+(\diver u)m=0,\\
			&m=(1-\alpha\Delta )u,\qquad u|_{t=0}=u_0,
		\end{aligned}\right.\qquad x\in\R^2.
	\end{equation}
	Here $u=(u_1,u_2)$ denotes the velocity, $m=(m_1,m_2)$ is the momentum, and the parameter $\alpha>0$ represents the square of a characteristic length scale. We adopt the convention $(\nabla u)_{ij}=\partial_j u_i$ and use the Einstein summation convention throughout. Following \cite{CL2012}, we introduce the stress tensor
	$$
	T_{ij}=m_iu_j+\frac 1 2\delta_{ij}|u|^2-\alpha\partial_i u\cdot\partial_j u+\frac 12 \alpha\delta_{ij}|\nabla u|^2$$
	and rewrite \eqref{EP} as 
	\begin{equation}\label{EP T}
		\partial_t m_i+\partial_j T_{ij}=0.
	\end{equation}
	For smooth solutions with enough spatial decay, there are two natural conservation laws
	\begin{equation}\label{conservation laws}
		\begin{array}{c}
			\frac{\dd}{{\dd}t}\int_{\R^2} m\,{\dd}x=0,\qquad
			\frac{\dd}{{\dd}t}\int_{\R^2}|u|^2+\alpha|\nabla u|^2\,{\dd}x=0.
		\end{array}
	\end{equation}
	
	\subsection{Background}
	
	The Euler--Poincar\'e equations were first introduced by Holm, Marsden, and Ratiu in \cite{HMR1998,HMR2002} and studied in \cite{HMR1998,HMR2002,HS2003} as a higher dimensional Camassa--Holm system for modeling and analyzing nonlinear shallow water waves. 
	
	In one dimension, the Euler--Poincar\'e equations reduce to the dispersionless Camassa--Holm (CH) equation
	$$
	\partial_t m+u\partial_xm+2m\partial_xu=0,\qquad m=(1-\alpha\partial_x^2)u,\qquad \alpha>0.
	$$
	Originally derived as a model for shallow water waves \cite{CH,CH2}, the CH equation possesses a bi-Hamiltonian structure, infinitely many conservation laws, and a rich integrable structure \cite{CH,CH6,CH4,CH5,CH7}. A distinctive feature is the coexistence of peaked solitary waves and wave breaking, in which the solution remains bounded while its slope becomes unbounded in finite time \cite{CH8,CH12,CH15}. For $\alpha=1$, it admits peakon solutions $u(t,x)=ce^{-|x-ct|}$ \cite{CH}, whose orbital stability was established by Constantin and Strauss \cite{CH17}.
	
	The Cauchy problem for the CH equation has been extensively studied in both the real-line and periodic settings \cite{CH8,CH9,CH13}. On the real line, local well-posedness holds in $H^s(\mathbb R)$ for $s>3/2$ \cite{CH11}. Danchin developed a Besov-space theory for $B^s_{p,r}(\mathbb R)$ with $1\leq p<\infty$, $1\leq r<\infty$, and $s>\max\{3/2,1+1/p\}$ \cite{CH10}, and subsequently established local well-posedness in the critical space $B^{3/2}_{2,1}(\mathbb R)$ and ill-posedness in $B^{3/2}_{2,\infty}(\mathbb R)$ \cite{CH18}. For weak solutions, global existence has been established via compensated compactness and viscous approximation, and weak--strong uniqueness as well as large-time behavior have been obtained under suitable assumptions \cite{CH19,CH20,CH21,CH22}. Thus, the one-dimensional problem is well understood for both strong solutions in critical spaces and global weak solutions.
	
	In two dimensions, the Euler-Poincar\'e equation is known as the averaged template matching equation in the computer vision literature \cite{HMA2001,HM2005,HS2009}. The local well-posedness in Sobolev spaces was initiated by Chae and Liu in \cite{CL2012}:
	\begin{itemize}
		\item Let $\alpha>0$ and $u_0\in H^k(\R^d)$ with $k>d/2+3$, then there exists $T=T(\norm{u_0}_{H^k})>0$ and $u\in C([0,T),H^k(\R^d))$, such that $u$ is a unique classical solution to the Euler-Poincar\'e equations.
		\item Let $\alpha=0$ and $u_0\in H^k(\R^d)$ with $k>d/2+1$, then there exists $T=T(\norm{u_0}_{H^k})>0$ and $u\in C([0,T),H^k(\R^d))$, such that $u$ is a unique classical solution to the Euler-Poincar\'e equations. Moreover, if $k>d/2+2$, $u_0(x)=-u_0(-x)$ and $\diver u_0(0)<0$, then the corresponding classical solution blows up in finite time.
	\end{itemize}
	Moreover, blow-up criteria, a zero--$\alpha$ limit and Liouville type theorem were established in \cite{CL2012}. Recently, the local well-posedness of the higher-dimensional CH equations in the Besov spaces $C([0,T);B^s_{p,r}(\mathbb{R}^d))$ was proved in \cite{YY2015} for $s>\max\{3/2,1+d/p\}$, and for $s=1+d/p$ with $r=1$.
	
	On conditions for blow-up and global existence of the Euler--Poincar\'e equations, Li, Yu and Zhai \cite{LYZ2013} assumed that $m=\nabla\phi$ with $\phi$ a radial scalar function. For $\alpha=1$, a direct computation reduces the system \eqref{EP} to
	\begin{equation}\label{eq:LYZ}
		\partial_t\phi(t,r)
		=\frac12\phi^2(r)
		+\int_r^\infty \partial_s\phi(t,s)\bigl[(1-\Delta)^{-1}\phi(t,s)\bigr]\,\dd s
		-\left((1-\Delta)^{-1}\nabla\phi\cdot\nabla\phi\right)(t,r).
	\end{equation}
	Based on \eqref{eq:LYZ}, they obtained blow-up and global existence results under suitable sign and monotonicity conditions, relying on the local well-posedness theory in \cite{CL2012}:
	\begin{itemize}
		\item Let $d\ge 2$, assume $\phi_0\in H^k(\R^d)$ for some $k>d/2+2$, in addition, assume $\phi_0\in\dot{B}^0_{1,\infty}(\R^2)$ for $d=2$. Let $u_0=(1-\Delta)^{-1}\nabla\phi_0$. If $\phi_0(0)\ge0$ and $\phi_0\not\equiv0$, then the corresponding solution blows up in finite time.
		\item Let $d\ge 2$, assume $\phi_0\in H^k(\R^d)$ for some $k>d/2+2$, in addition, assume $\phi_0\in\dot{B}^0_{1,\infty}(\R^2)$ for $d=2$. Let $u_0=(1-\Delta)^{-1}\nabla\phi_0$. If $\phi_0(0)\le0$ and $\phi_0'(r)\ge 0$ for any $r\in[0,\infty)$, then the corresponding solution exists globally in time.
	\end{itemize}
	For the periodic Euler--Poincar\'e equations, Luo and Yin \cite{LY2020} established finite-time blow-up for a special class of smooth initial data satisfying directional invariance conditions. Their analysis exploits the rotational covariance of the system to obtain a reduced evolution structure. They also investigated ill-posedness in critical Besov spaces and constructed periodic peakon solutions.
	
	The radial gradient-momentum ansatz in \cite{LYZ2013} yields a monotonicity identity for the central potential, but it excludes a nonzero tangential momentum component, i.e., a nonzero solenoidal component. This raises the question of whether the same blow-up mechanism persists when such a component is present. In this paper, we address this question on $\mathbb{R}^2$ based on a div--curl decomposition with radial scalar potentials. This structure preserves rotational equivariance while allowing nonzero tangential momentum. The resulting angular coupling introduces a nonpositive contribution to the evolution of the central potential, whose control constitutes the main additional difficulty in our analysis.
	
	\subsection{Main results}
	
	Inspired by the special structure introduced in \cite{CL2012}, we consider the following div--curl decomposition:
	\begin{equation}
		m=\nabla\phi_1+\nabla^\perp\phi_2,
	\end{equation}
	with $\nabla^\perp=(-\partial_2,\partial_1).$  Here~$\phi_1$ and $\phi_2$ are radial scalar-valued functions. In particular,
	\begin{equation}\label{divcurl}
		\diver m=\Delta\phi_1,\qquad \curl m=\Delta\phi_2.
	\end{equation}
	
	Introducing polar coordinates $x=(r\cos\theta,r\sin\theta)$, where $r=|x|$ and $\theta\in[0,2\pi)$, we define the radial and angular unit vectors by
	$$
	e_1=(\cos\theta,\sin\theta),
	\qquad
	e_2=(-\sin\theta,\cos\theta).
	$$
	For radial scalar potentials $\phi_1$ and $\phi_2$, we then have
	$$
	\nabla\phi_1=(\partial_r\phi_1)e_1,
	\qquad
	\nabla^\perp\phi_2=(\partial_r\phi_2)e_2.
	$$
	Set $p(t,r)=\partial_r\phi_1(t,r)$ and $q(t,r)=\partial_r\phi_2(t,r)$. With the normalization $\phi_i(t,r)\to0$ as $r\to\infty$, the potentials admit the representations
	\begin{equation}\label{potentials}
		\phi_1(t,r)=-\int_r^\infty p(t,s)\dd s,\qquad
		\phi_2(t,r)=-\int_r^\infty q(t,s)\dd s.
	\end{equation}
	
	Fix $r_0>0$ and $k>3$, and consider the initial radial potentials satisfying
	\begin{equation}\label{data}
		\begin{aligned}
			&\phi_{1,0}\in H^k(\R^2),\qquad q_0\in C_0^\infty([r_0,\infty)),\\
			&\phi_{2,0}(r)=-\int_r^\infty q_0(s)\dd s,\qquad
			q_0=\frac{\dd}{\dd r}\phi_{2,0}(r).
		\end{aligned}
	\end{equation}
	and set
	\begin{equation}\label{u0}
		m_0=\nabla\phi_{1,0}+\nabla^\perp\phi_{2,0},\qquad
		u_0=(1-\Delta)^{-1}m_0,\qquad E=\norm{u_0}_{H^1}.
	\end{equation}
	
	For simplicity, we set $\alpha=1$ throughout the remainder of the paper. We establish a sufficient condition for finite-time blow-up and prove ill-posedness in the critical Sobolev space $H^2(\mathbb R^2)$. Our main results are stated as follows:
	
	\begin{theo}\label{thm:blow}
		Let $\phi_{1,0}$ and $\phi_{2,0}$ be radial functions satisfying \eqref{data} for a fixed $r_0>0$, and let $m_0$ and $u_0$ be defined by \eqref{u0}. Denote by $T^*$ the maximal existence time of the corresponding solution to \eqref{EP}. There exists a constant $C_0=C_0(r_0)>0$ such that, if
		\begin{equation}\label{criterion}
			A_0:=\phi_{1,0}(0)>
			C_0\bigl(\norm{\phi_{1,0}}_{L^2}+E\bigr),
		\end{equation}
		then
		\begin{equation}\label{lifespan}
			T^*\leq\frac{8}{A_0}<\infty.
		\end{equation}
	\end{theo}
	
	\begin{theo}[Norm inflation in $H^2$]\label{thm:inflation}
		For any $0<\eps<1$, there exists $u_0\in H^\infty(\R^2)$ such that the following holds:
		\begin{enumerate}
			\item $\|u_0\|_{H^2(\R^2)}\le \varepsilon$;
			\item There exists a time $0<t_\varepsilon<\min\{\eps,T^*\}$ such that the corresponding solution $u\in C([0,T^*);H^\infty(\R^2))$ satisfies
			\begin{equation}\label{inflation}
				\norm{u(t_\eps)}_{H^2}>\eps^{-1}.
			\end{equation}
		\end{enumerate}
	\end{theo}
	
	Compared with the two-dimensional radial gradient-momentum setting in \cite{LYZ2013}, we consider a div--curl decomposition of the momentum with radial scalar potentials, allowing nonzero tangential momentum while preserving rotational equivariance. For this class of solutions, we establish a quantitative sufficient condition for finite-time blow-up: the initial central potential $\phi_{1,0}(0)$ is required to dominate the initial $L^2$ norm of the irrotational potential and the conserved $H^1$ norm of the velocity. Under this condition, we obtain an explicit upper bound for the maximal existence time in terms of $\phi_{1,0}(0)$; see Theorem~\ref{thm:blow}. Furthermore, we establish norm inflation at the zero initial datum in the critical Sobolev space $H^2(\mathbb R^2)$: arbitrarily small smooth initial data in $H^2$ can generate solutions with arbitrarily large $H^2$ norms in arbitrarily short times; see Theorem~\ref{thm:inflation}.
	
	The main idea in proving Theorem~\ref{thm:blow} is to control the non-positive contribution arising from the tangential momentum component $\nabla^\perp\phi_2=q e_2$ in the evolution equation for the central potential $A(t):=\phi_1(t,0)$. To overcome this difficulty, we use the support propagation property of the homogeneous transport equation satisfied by the tangential momentum $q$ to show that it remains zero near the origin on a short time interval determined by the conserved energy and the inner radius of the support of $q_0$. Moreover, the corresponding tangential velocity $(1-\Delta)^{-1}q e_2$ allows us to apply interior elliptic estimates to the tangential velocity component in the region where $q$ vanishes. Combining these estimates with the energy bound in the exterior region yields control of the singularly weighted tangential contribution. Once this estimate is obtained, since the divergence-free component disappears after applying the gradient projection, we use the conservation-form momentum equation \eqref{EP T} to control $\phi_1$ by the $\norm{u}_{H^1}$ and $\norm{\phi_{1,0}}_{L^2}$. These estimates yield a Riccati-type inequality for the central potential $\phi_1(t,0)$, implying finite-time breakdown of the classical solution.
	
	Regarding critical ill-posedness, we construct radial potentials with disjoint frequency supports, so that the squared $L^2$ norms of their gradients add while their values at the origin add with the same sign. This allows the initial velocity to have an arbitrarily small $H^2$ norm while the initial central potential is arbitrarily large. The quantitative blow-up criterion then yields an arbitrarily short maximal existence time. Combining the deformation-tensor blow-up criterion in \cite{CL2012} with the estimate $\|S\|_{\dot B^0_{\infty,\infty}}\lesssim\|u\|_{H^2}$, we conclude that the $H^2$ norm can be made arbitrarily large before the maximal existence time, which yields norm inflation near the zero initial datum.
	
	\section{Notations and preliminaries}\label{sec:prelim}
	In this section, we introduce some notations and recall the local well-posedness result for the Euler--Poincar\'e equations. 
	
	We first recall the Littlewood--Paley decomposition \cite{BCD}. Choose smooth cutoffs $\chi,\varphi$ with values in $[0,1]$, with $\chi$ supported in $|\xi|\leq 4/3$ and $\varphi$ supported in $3/4\leq|\xi|\leq 8/3$. They can be chosen so that
	$$
	\forall\xi\in\mathbb{R}^d, \chi(\xi)+\sum_{j\geq0}\varphi(2^{-j}\xi)=1,
	$$
	and
	$$
	{\rm supp}\ \varphi(2^{-j}\cdot)\cap {\rm supp}\ \varphi(2^{-j'}\cdot)=\emptyset,\ {\rm if}\  |j-j^{\prime}|\geq2,
	$$
	$$
	{\rm supp}\ \chi(\cdot)\cap {\rm supp}\ \varphi(2^{-j}\cdot)=\emptyset,\quad \text{if }j\geq1.
	$$
	Then we define the Littlewood-Paley operators
	\begin{equation}\label{Littlewood Paley}
		\widehat{P_{j}f}(\xi):=\varphi(2^{-j}\xi)\hat{f}(\xi),\quad
		\widehat{P_{>j}f}(\xi):=(1-\chi(2^{-j}\xi))\hat{f}(\xi),\quad
		\widehat{P_{\le j}f}(\xi):=\chi(2^{-j}\xi)\hat{f}(\xi).
	\end{equation}
	We need the Bernstein inequality. For any $1\le p<q\le\infty$,
	$$
	\begin{array}{c}
		\norm{P_{j}f}_{L^q(\R^d)}\le C2^{jd(1/p-1/q)}\norm{P_jf}_{L^p(\R^d)},\\
		\norm{P_{\le j}f}_{L^q(\R^d)}\le C2^{jd(1/p-1/q)}\norm{P_{\le j}f}_{L^p(\R^d)}.
	\end{array}
	$$
	\subsection{Local well-posedness and blow-up criteria}
	In order to show that the solutions satisfy the rotational invariance, we need  the local well-posedness result:
	\begin{lemm}\cite{CL2012}
		Let $u_0\in H^k(\R^d)$ with $k>\frac d 2+3$. Then there exists $T=T(\norm{u_0}_{H^k(\R^d)})$ such that the Euler--Poincar\'e equations have a unique solution $u\in C([0,T];H^k(\R^d))$.
	\end{lemm}

	Therefore, the uniqueness preserves rotational invariance. Consequently, according to \eqref{potentials}, we assume that the velocity $u$ and the momentum $m$ admit the following radial--tangential decompositions:
	\begin{equation}\label{polar}
		u=v(t,r)e_1+w(t,r)e_2,\qquad
		m=p(t,r)e_1+q(t,r)e_2.
	\end{equation}
	According to $\eqref{EP}_2$, a direct calculation yields
	\begin{equation}\label{elliptic}
		\begin{array}{c}
			p=\A v,\quad q=\A w,\quad\text{with}\quad \A=1-\partial_r^2-r^{-1}\partial_r+r^{-2}.\\
			v=(1-\Delta)^{-1}\partial_r\phi_1,\quad w=(1-\Delta)^{-1}\partial_r\phi_2.
		\end{array}
	\end{equation}
	
	In the orthonormal polar coordinate,
	$$
	\nabla u=\begin{pmatrix}v_r&-w/r\\w_r&v/r\end{pmatrix}.
	$$
	Substitution into \eqref{EP} yields
	\begin{equation}\label{p q-eq}
		\left\{
		\begin{aligned}
			&p_t+vp_r+(2v_r+v/r)p+(w_r-w/r)q=0,\\
			&q_t+vq_r+(v_r+2v/r)q=0.
		\end{aligned}
		\right.
	\end{equation}
	Assuming $\phi_{i,0}$ decays sufficiently fast at infinity, we may integrate \eqref{p q-eq} and obtain
	\begin{equation}\label{phi-eq}
		\left\{\ 
		\begin{aligned}
			\partial_t\phi_1(t,r)=&\int_r^\infty\partial^2_s\phi_1(t,s)(1-\Delta)^{-1}\partial_s\phi_1(t,s)\dd s\\
			&+\int_r^\infty\left(2(1-\Delta)^{-1}\partial^2_s\phi_1+\frac{(1-\Delta)^{-1}\partial_s\phi_1}{s}\right)\partial_s\phi_1(t,s)\dd s\\
			&+\int_r^\infty\left((1-\Delta)^{-1}\partial_s^2\phi_2-\frac{(1-\Delta)^{-1}\partial_s\phi_2}s\right)\partial_s\phi_2(t,s)\dd s,\\
			\partial_t\phi_2(t,r)=&\int_r^\infty\partial_s^2\phi_2(t,s)(1-\Delta)^{-1}\partial_s\phi_1(t,s)\dd s\\
			&+\int_r^\infty\left((1-\Delta)^{-1}\partial_s^2\phi_1+\frac{2(1-\Delta)^{-1}\partial_s\phi_1}{s}\right)\partial_s\phi_2(t,s)\dd s.
		\end{aligned}
		\right.
	\end{equation}
	In particular, the tangential momentum coefficient $q$ satisfies a homogeneous transport equation. Besides, the conservation law $\eqref{conservation laws}_2$ can be rewritten as
	\begin{equation}\label{energy}
		\norm{u(t)}_{H^1}^2=E^2
		=2\pi\int_0^\infty\left(v^2+w^2+v_r^2+w_r^2+
		\frac{v^2+w^2}{r^2}\right)r\dd r.
	\end{equation}
	
	\begin{prop}\label{prop:well posedness}
		Assume initially $m_0=\nabla\phi_{1,0}+\nabla^\perp\phi_{2,0}$, where $\phi_{i,0}\in H^k$ are radial functions on $\R^2$ for some $k>3$. Then for any $t>0$, the solution $m(t)=(1-\Delta)u(t)$ can be written as $m(t)=\nabla\phi_{1}(t)+\nabla^\perp\phi_{2}(t)$. Then $\phi_{1}(t,r)$ and $\phi_{2}(t,r)$ solve \eqref{phi-eq} in the classical sense. Moreover, $\phi_1$ satisfies the growth estimate
		\begin{equation}\label{phi1 growth}
			\norm{\phi_1(t)}_{L^2(\R^2)}\le C_1(1+t), \quad\forall t\ge0.
		\end{equation}
		Here, $C_1$ depends only on initial data $\phi_{i,0}$ for $i=1,2$.
	\end{prop}
	\begin{proof}
		Since $m=\nabla\phi_{1}+\nabla^\perp\phi_{2}$, we have $u=(1-\Delta)^{-1}\nabla\phi_{1}+(1-\Delta)^{-1}\nabla^\perp\phi_{2}$. 
		According to the div--curl decomposition and the conservation law of $u$,
		$$\norm{u}_{L^2}^2=\norm{(1-\Delta)^{-1}\nabla\phi_1}_{L^2}^2+\norm{(1-\Delta)^{-1}\nabla^\perp\phi_2}_{L^2}^2\le E^2,$$
		thus, by the Bernstein inequality
		$$\norm{P_{\ge 1}\phi_1}_{L^2}\le C\norm{P_{\ge 1}(1-\Delta)^{-1}\nabla\phi_1}_{H^1}\le CE.$$
		
		For the low frequency part, by \eqref{EP T} and $\eqref{divcurl}_1$, we obtain that
		$$\phi_1(t,r)=\phi_{1,0}(r)-\int_0^t\Delta^{-1}\partial_i\partial_jT_{ij}(\tau)\dd\tau,$$
		therefore, by the Bernstein inequality, we have
		$$\begin{aligned}
			\norm{P_{<1}\phi_{1}}_{L^2}
			&\le
			\norm{P_{<1}\phi_{1,0}}_{L^2}+\int_0^t\norm{P_{<1}\Delta^{-1}\partial_i\partial_jT_{ij}}_{L^2}(\tau)\dd \tau\\
			&\le
			\norm{\phi_{1,0}}_{L^2}+C\int_0^t\norm{u}_{H^1}^2(\tau)\dd \tau
			\le 
			\norm{\phi_{1,0}}_{L^2}+CE^2t.
		\end{aligned}
		$$
		
		Combining the low and high frequency estimates, we obtain \eqref{phi1 growth}.
	\end{proof}

	The following blow-up criteria is used to prove the blow up and ill-posedness.
	\begin{lemm}\label{lem:blow up criteria}\cite{CL2012}
		Let $u_0\in H^k(\R^d)$ with $k>\frac d 2+3$. Suppose that $T^*$ is the lifespan of the solution $u$ to the Euler-Poincar\'e equations. If $T^*<\infty$, then
		$$\int_0^{T^*}\norm{S(t)}_{\dot{B}^0_{\infty,\infty}(\R^d)}\,\dd t=\infty.$$
	\end{lemm}
	
	Analogously to \cite[Lemma~1.2]{LYZ2013}, based on \eqref{p q-eq}, we establish the following blow-up criterion for the two-dimensional Euler--Poincar\'e equations under the div--curl structure:
	\begin{lemm}\label{lem:continuation}
		Let $\phi_{1,0},\phi_{2,0}\in H^k(\R^2)$ be radial for some $k>3$. Let $u$ be the maximal lifespan solution corresponding to initial data $u_0=(1-\Delta)^{-1}m_0=(1-\Delta)^{-1}(\nabla\phi_{1,0}+\nabla^\perp\phi_{2,0})$. If the maximal lifespan $T^*$ is finite, then
		$$
		\int_0^{T^*}\norm{\phi_1(t)}_{L^\infty}\dd t=\infty.
		$$
	\end{lemm}
	
	\begin{proof}
		Assume that $T^*<\infty$ and $\int_0^{T^*}\norm{\phi_1(t)}_{L^\infty}\dd t<\infty$. According to \eqref{polar}, we can easily obtain that $u(t,-x)=-u(t,x)$, due to the continuity of $u$ at $x=0$, $v(t,0)=w(t,0)=0$.
		
		Set $g=(1-\Delta)^{-1}\phi_1$. Since $\partial_1^2+\partial_2^2=\partial_r^2+r^{-1}\partial_r+r^{-2}\partial_\theta^2$, we deduce that $\phi_1=g-g''-r^{-1}g'$. Since $v=\partial_r g$, we have
		\begin{equation}\label{g-phi1}
			v_r+\frac{v}{r}=g-\phi_1.
		\end{equation}
		Using the regularity at the origin, we obtain
		$$
		\frac{v(t,r)}{r}
		=
		\frac{1}{r^2}\int_0^r
		s\bigl(g(t,s)-\phi_1(t,s)\bigr)\dd s.
		$$
		The $L^\infty$ boundedness of $(1-\Delta)^{-1}$ gives $\norm{g(t)}_{L^\infty}
		\leq\norm{\phi_1(t)}_{L^\infty}$.
		Consequently,
		$$
		\norm{v(t)/r}_{L^\infty}
		\leq\norm{\phi_1(t)}_{L^\infty},
		\qquad
		\norm{v_r(t)}_{L^\infty}
		\leq C\norm{\phi_1(t)}_{L^\infty}.
		$$
		
		Thus, the $L^\infty$-estimate of the transport equation $\eqref{p q-eq}_2$ yields
		\begin{equation}\label{q infty}
			\norm{q(t)}_{L^\infty}
			\leq
			\norm{q_0}_{L^\infty}
			\exp\left(
			\int_0^t\norm{v_r+2v/r}_{L^\infty}\dd s
			\right)
			\le
			\norm{q_0}_{L^\infty}
			\exp\left(
			C\int_0^t\norm{\phi_1(s)}_{L^\infty}\dd s
			\right).
		\end{equation}
		
		Let $G$ denote the scalar convolution kernel of $(1-\Delta)^{-1}$. Since $\nabla G\in L^1(\R^2)$, the Young inequality implies
		$$
		\norm{\nabla(w e_2)}_{L^\infty}=\norm{\nabla G*(q e_2)}_{L^\infty}
		\leq
		\norm{\nabla G}_{L^1}\norm{q e_2}_{L^\infty}
		\leq C\norm{q}_{L^\infty}.
		$$
		On the other hand, the radial velocity satisfies
		$$
		\norm{\nabla(v e_1)}_{L^\infty}
		\leq
		C\left(
		\norm{v_r}_{L^\infty}
		+\norm{v/r}_{L^\infty}
		\right)\le C\norm{\phi_1(t)}_{L^\infty}.
		$$
		Combining the above two inequalities and \eqref{q infty}, we obtain
		$$
		\norm{\nabla u(t)}_{L^\infty}
		\le
		\norm{\nabla(w e_2)}_{L^\infty}+\norm{\nabla(v e_1)}_{L^\infty}
		\leq
		C\left(
		\norm{\phi_1(t)}_{L^\infty}
		+\norm{q_0}_{L^\infty}\exp\left(
		C\int_0^t\norm{\phi_1(s)}_{L^\infty}\dd s
		\right)
		\right).
		$$
		Therefore,
		$$
		\int_0^{T^*}\norm{\nabla u(t)}_{L^\infty}\dd t
		\leq
		C\int_0^{T^*}\norm{\phi_1(t)}_{L^\infty}\dd t+
		CT^*\norm{q_0}_{L^\infty}\exp(C\int_0^{T^*}\norm{\phi_1(t)}_{L^\infty}\dd t)
		<\infty.
		$$
		
		Finally, let $S=(\nabla u+(\nabla u)^T)/2$. The embedding
		$L^\infty\hookrightarrow\dot B^0_{\infty,\infty}$ gives
		$$
		\int_0^{T^*}\norm{S(t)}_{\dot B^0_{\infty,\infty}}\dd t
		\leq
		C\int_0^{T^*}\norm{\nabla u(t)}_{L^\infty}\dd t
		<\infty.
		$$
		This contradicts Lemma \ref{lem:blow up criteria} and completes the proof.
	\end{proof}
	\begin{rema}
		When $q=0$, the tangential momentum vanishes and
		$m=\nabla\phi_1$. Thus, Lemma~\ref{lem:continuation} reduces to
		the blow-up criterion for purely radial solutions in
		\cite[Lemma~1.2]{LYZ2013}.
	\end{rema}
	
	\subsection{Finite propagation speed}
	In order to obtain the blow-up result, the main difficulty comes from the quantity $(w_r-w/r)q$ in \eqref{p q-eq}. To overcome this, we need to use the finite propagation speed property.
	In the following lemma, we use the conserved energy to obtain an explicit time $T_{r_0}$ such that $q(t,r)=0$ for $0<r<\frac{r_0}{2}$.
	
	\begin{lemm}\label{lem:core}
		Let $q_0\in C_0^\infty([r_0,\infty))$. Assume that $\phi_{1,0}$ and $\phi_{2,0}$ are radial functions satisfying \eqref{data} for a fixed $r_0>0$, and let $m_0$ and $u_0$ be defined by \eqref{u0}. Denote by $T^*$ the maximal $H^{k+1}$ existence time of the corresponding solution to \eqref{EP}. Then there exists
		\begin{equation}\label{protectiontime}
			T_{r_0}=
			\begin{cases}
				\dfrac{\sqrt{2\pi}\,r_0}{2E},&E>0,\\[4pt]
				+\infty,&E=0,
			\end{cases}
		\end{equation}
		such that
		\begin{equation}\label{core}
			q(t,r)=0\quad(0<r<r_0/2),\qquad
			0\leq t<\min\{T^*,T_{r_0}\}.
		\end{equation}
	\end{lemm}
	\begin{proof}
		If $E=0$, the solution is identically zero, and the conclusion is immediate. For $E>0$, consider the radial characteristic
		$$R_t(t,a)=v(t,R(t,a)),\quad R(0,a)=a,$$
		it is easy to deduce that
		$$
		R_a(t,a)=\exp\left(\int_0^t v_r(s,R(s,a))\dd s\right)>0.
		$$
		According to $\eqref{p q-eq}_2$, we have
		$$
		\frac{\dd}{\dd t}(R^2(t,a)R_a(t,a)q(t,R(t,a)))=0.
		$$
		Thus, combining this identity with $v(t,0)=0$ we deduce that for $a>0$,
		\begin{equation}\label{qchar}
			q(t,R(t,a))R_a(t,a)\left(\frac{R(t,a)}a\right)^2=q_0(a).
		\end{equation}
		
		By \eqref{energy},
		$$
		\int_0^\infty\left(v_r^2+\frac{v^2}{r^2}\right)r\dd r
		\leq\frac{E^2}{2\pi},
		$$
		the fundamental theorem of calculus implies that for all $r>0$,
		\begin{align*}
			|v(t,r)|^2
			&=2\int_0^r v(t,r')v_r(t,r')\dd r'\\
			&\leq\int_0^r\left(\frac{|v(t,r')|^2}{r'}
			+r'|v_r(t,r')|^2\right)\dd r'
			\leq\frac{E^2}{2\pi}.
		\end{align*}
		Therefore
		\begin{equation}\label{radial-speed}
			\norm{v(t)}_{L^\infty(0,\infty)}
			\leq \frac{E}{\sqrt{2\pi}}.
		\end{equation}
		Consequently every radial characteristic satisfies
		\begin{equation}\label{support-cone}
			|R(t,a)-a|\leq \frac{E}{\sqrt{2\pi}}t,\qquad R(t,a)\geq a-\frac{E}{\sqrt{2\pi}}t.
		\end{equation}
		Then if $q(t,R(t,a))\ne0$, \eqref{qchar} implies $q_0(a)\ne0$, hence $a\geq r_0$. Thus, for $0\leq t<\min\{T^*,T_{r_0}\}$,
		$$
		q(t,r)=0,\qquad 0<r<r_0/2.
		$$
		This completes the proof.
	\end{proof}

	\section{Finite-time blow-up}
	In this section we prove the main theorem.
	\begin{proof}[Proof of Theorem \ref{thm:blow}]
		According to Proposition \ref{prop:well posedness}, $\phi_1(t),\phi_2(t)\in H^k$, then we can deduce that
		\begin{equation}\label{v w r=0}
			v(t,0)=w(t,0)=0,\quad \frac{v(t,r)}{r}\to v_r(t,0)\quad\text{and}\quad \frac{w(t,r)}{r}\to w_r(t,0)\quad\text{as}\ r\to0^{+}.
		\end{equation}
		Combining \eqref{v w r=0} with \eqref{g-phi1}, $g(t,0)-\phi_1(t,0)=2v_r(t,0)$. By \eqref{p q-eq}:
		$$
		\frac{\dd}{\dd t}\phi_1(t,0)
		=
		-\int_0^\infty \frac{\partial}{\partial t}p(t,r)\,\dd r\\
		=
		\int_0^\infty \left(v_r+\frac v r\right)p(t,r)+(w_r-\frac w r)q(t,r)\dd r.
		$$
		According to \eqref{elliptic} and \eqref{v w r=0}, $p=\A v=v-\partial_r(v_r+\frac v r)$ and $q=\A w=w-\partial_r\left(w_r-\frac w r\right)-\frac 2 r\left(w_r-\frac w r\right)$, integrating by parts, we have
		$$
		\begin{aligned}
			\int_0^\infty\left(v_r+\frac v r\right)p(t,r)\dd r
			=&
			\int_0^\infty \left(v_r+\frac v r\right) v\dd r-\int_0^\infty \left(v_r+\frac v r\right)\partial_r\left(v_r+\frac v r\right)\dd r\\
			=&
			\int_0^\infty \frac{v^2} r\dd r -\frac12\left(v_r+\frac vr\right)^2\bigg|_0^\infty
			=
			\int_0^\infty \frac{v^2} r\dd r+\frac 1 2(\phi_1(t,0)-g(t,0))^2,
		\end{aligned}
		$$
		and 
		$$
		\begin{aligned}
			\int_0^\infty(w_r-\frac w r)q(t,r)\dd r
			=&
			-2\int_0^\infty\frac 1 r\left(w_r-\frac w r\right)^2\dd r+\int_0^\infty w\left(w_r-\frac w r\right)\dd r -\frac12\left(w_r-\frac wr\right)^2\bigg|_0^\infty\\
			=& -\int_0^\infty\frac{w^2}{r}+\frac2r\left(w_r-\frac wr\right)^2\dd r.
		\end{aligned}
		$$
		Then we have
		\begin{equation}\label{phi_1(t,0)-eq}
			\begin{aligned}
				\frac{\dd}{\dd t}\phi_1(t,0)
				&=
				\frac 1 2(\phi_1(t,0)-g(t,0))^2	+\int_0^\infty\frac {v^2} r\dd r -\int_0^\infty\frac{w^2}{r}+\frac2r\left(w_r-\frac wr\right)^2\dd r\\
				&=:
				\frac 1 2(\phi_1(t,0)-g(t,0))^2	+\int_0^\infty\frac {v^2} r\dd r-I.
			\end{aligned}
		\end{equation}
		
		According to conservation law \eqref{energy}, 
		\begin{equation}\label{I high}
			\int_{r_0/4}^\infty \frac 2 r\left(w_r-\frac w r\right)^2\dd r\le C_{r_0}\int_{r_0/4}^\infty \frac{w^2} r+rw_r^2\dd r\le C_{r_0}E^2.
		\end{equation}
		Since $q_0(r)=0$ for $r<r_0$, Lemma~\ref{lem:core} ensures that there exists a time $T_{r_0}>0$ such that
		$$
		q(t,r)=0
		\qquad\text{for }0<r<\frac{r_0}{2}
		\quad\text{and}\quad
		0\le t<\min\{T^*,T_{r_0}\}.
		$$
		Therefore, $q e_2=(1-\Delta)(we_2)=0$ in the ball $B_{r_0/2}$ for $0\le t<\min\{T^*,T_{r_0}\}$. By the interior $H^4$ estimate for uniformly elliptic equations \cite{Evans2010},
		applied componentwise to $(1-\Delta)(we_2)=0$
		in $B_{r_0/2}$ and the Sobolev embedding
		$H^4(\R^2)\hookrightarrow C^2(\R^2)$, we obtain
		$$
		\|we_2\|_{C^2(\overline{B_{r_0/4}})}
		\leq C_{r_0}\|we_2\|_{H^4(B_{r_0/4})}
		\leq C_{r_0}\|we_2\|_{L^2(B_{r_0/2})}
		\leq C_{r_0}E.
		$$
		Since $w(0)=0$, the fundamental theorem of calculus yields
		$$
		\begin{aligned}
			\left|w_r(r)-\frac{w(r)}r\right|
			&\leq \frac1r\int_0^r|w_r(r)-w_r(s)|\,ds\\
			&\leq
			\frac{\|w_{rr}\|_{L^\infty(B_{r_0/4})}}r
			\int_0^r(r-s)\,ds
			\leq C_{r_0}Er,
		\end{aligned}
		$$
		for $0<r<r_0/4$. Thus, we have
		\begin{equation}\label{I low}
			\int_0^{r_0/4}\frac 2 r\left(w_r-\frac w r\right)^2\dd r\le C_{r_0}\int_0^{r_0/4} E^2r\dd r\le C_{r_0}E^2.
		\end{equation}
		Combining \eqref{I high}, \eqref{I low} and \eqref{energy}, we deduce that
		\begin{equation}\label{I estiamte}
			I\le C_{r_0}E^2.
		\end{equation}
		
		Hence, substituting \eqref{I estiamte} into \eqref{phi_1(t,0)-eq}, we obtain that for all $0\le t<\min\{T^*,T_{r_0}\}$
		\begin{equation}\label{phi1 inequality}
			\frac{\dd}{\dd t}\phi_{1}(t,0)\ge \frac 1 2(\phi_1(t,0)-g(t,0))^2-C_{r_0}E^2.
		\end{equation}
		
		Denote the Bessel potential $K(x)=\mathcal F^{-1}\left((1+|\xi|^2)^{-1}\right)(x)$, and $K\in L^1(\R^2)\cap L^q(\R^2)$ for any $q<\infty$. According to Proposition \ref{prop:well posedness}, for $0\le t<\min\{T^*,T_{r_0}\}$,
		$$|g(t,0)|\le \norm{K}_{L^2}\norm{\phi_1}_{L^2}\le C_1\norm{K}_{L^2}(1+t)\le C_g.$$
		Choose $C_0$ sufficiently large such that
		$$C_0>
		\max\left\{
		\frac{4C_g}{\norm{\phi_{1,0}}_{L^2}+E},
		4C_{r_0}^{1/2},
		\frac 8 {ET_{r_0}}
		\right\},
		$$
		then we can deduce that
		\begin{equation}\label{C0 bound}
			C_{r_0}E^2\le \frac{A_0^2}{16},\quad |g(t,0)|\le \frac{A_0} 4,\quad \frac 8{A_0}<T_{r_0}.
		\end{equation}
		
		Denote $A(t)=\phi_1(t,0)$, since $A(0)=\phi_{1,0}(0)=A_0>0$, substituting \eqref{C0 bound} into \eqref{phi1 inequality} yields $A'(0)>0$. Similarly, a standard continuity argument therefore gives
		$$
		A(t)\geq A_0,\quad
		A'(t)\geq\frac 1 8A(t)^2,
		\quad
		0\leq t<\min\{T^*,T_{r_0}\}.
		$$ 
		Thus,
		$$
		\frac1{A(t)}\leq\frac1{A_0}-\frac t8,
		\qquad 0\leq t<\min\{T^*,T_{r_0}\}.
		$$
		Since $8/A_0<T_{r_0}$, the above inequality forces
		$\phi_1(t,0)$ to blow up as $t$ tends to $8/A_0$. Hence $T^*\leq8/A_0$.
	\end{proof}
	
	\begin{rema}\label{rem:ARMA}
		If $q_0=0$, the homogeneous transport equation gives $q(t)=0$, hence $w(t)=0$ and $I=0$ throughout the classical lifespan. In this case the concentration assumption \eqref{criterion} can be replaced by
		\begin{equation}\label{ARMA-data}
			\phi_{1,0}\in H^k(\R^2)\cap\dot B^0_{1,\infty}(\R^2),\qquad
			k>3,\qquad \phi_{1,0}\not\equiv0,\qquad \phi_{1,0}(0)\geq0 .
		\end{equation}
		The low-frequency hypothesis in \eqref{ARMA-data} yields the sharper growth bound
		$\norm{\phi_1(t)}_{L^2(\R^2)}\leq C_0(1+t)^{1/2}$ from \cite[Proposition~1.1]{LYZ2013}.
		
		The reduction is also visible from our estimates. The angular contribution must be retained as exactly zero, rather than replaced by the coarse error $Ce^2$. Formula \eqref{phi_1(t,0)-eq} then becomes
		$$
		\frac{\dd}{\dd t}\phi_1(t,0)=\int_0^\infty\frac{|g_r|^2}{r}\dd r
		+\frac12(\phi_1(t,0)-g(t,0))^2\geq0 .
		$$
		Thus, when $q_0=0$, we recover the blow-up criterion in \cite{LYZ2013}
		$$
		\phi_{1,0}\not\equiv0,
		\qquad
		\phi_{1,0}(0)\geq0.
		$$
		However, for $q_0\not\equiv0$, our argument still requires
		the stronger condition \eqref{criterion}.
	\end{rema}
	
	\section{Ill-posedness in critical Sobolev spaces}\label{sec:ill}
	\begin{proof}[Proof of Theorem \ref{thm:inflation}]
		Choose a real radial Schwartz function $\rho$ whose Fourier transform is nonnegative, nonzero, and supported in $1\leq|\xi|\leq2$. Let $\rho(0)>0$ and fix a nonzero real $\chi\in C_c^\infty((r_0,2r_0))$. For $\delta>0$ and $N\geq1$, define
		\begin{equation}\label{phi 0}
			\phi_{1,0}^N(x)=\frac{\delta}{\sqrt N}\sum_{j=1}^N\rho(2^{4j}x),
			\quad
			\phi_{2,0}^N(r)=-\delta\int_r^\infty\chi(s)\,\dd s.
		\end{equation}
		Let $u_0^N$ be given by \eqref{u0}, and the initial data satisfy \eqref{data}. Since the frequency support of $\rho$ are disjoint,
		\begin{equation}\label{phi10 initial}
			\begin{aligned}
				&\phi_{1,0}^N(0)=\delta\sqrt{N}\,\rho(0),\quad
				\|\nabla\phi_{1,0}^N\|_{L^2}^2=\delta^2\|\nabla\rho\|_{L^2}^2,\\
				&\quad\ \|\phi_{1,0}^N\|_{L^2}^2=\frac{\delta^2}{N}\|\rho\|_{L^2}^2
				\sum_{\ell=1}^N2^{-2j_\ell}\leq C\delta^2/N.
			\end{aligned}
		\end{equation}
		It is obvious that $\norm{\phi_{2,0}}_{H^2}\le C\delta$. Since $u_0^N=(1-\Delta)^{-1}\nabla\phi_{1,0}^N+(1-\Delta)^{-1}\nabla^\perp\phi_{2,0}^N$. By Plancherel's theorem and the Fourier multiplier
		representation of $(1-\Delta)^{-1}$,
		$$
		\|(1-\Delta)^{-1}\nabla f\|_{H^2}=\|\nabla f\|_{L^2},
		$$
		then we have
		\begin{equation}\label{u0 H^2}
			\norm{u_0^N}_{H^2}\le \norm{\nabla\phi_{1,0}^N}_{L^2}+\norm{\nabla\phi_{2,0}^N}_{L^2}\le C_1\delta.
		\end{equation}
		Besides, by direct computation
		\begin{equation}\label{u0 H^1}
			\norm{u_0^N}_{H^1}\le C\left(\norm{\phi_{1,0}^N}_{L^2}+\norm{\phi_{2,0}}_{L^2}\right)\le C_2 \delta(1+N^{-1/2}).
		\end{equation}
		The constants in these estimates are independent of $\delta$ and $N$.
		
		Fix $0<\eps<1$ and choose $\delta>0$ with $C_1\delta<\eps$. For this fixed $\delta$, choose a finite $N$ sufficiently large such that
		$$
		\phi_{1,0}^N(0)>C_0(r_0)(\|\phi_{1,0}^N\|_{L^2}+\|u_0^N\|_{H^1}),
		\qquad \frac{8}{\delta\sqrt{N}\rho(0)}<\eps.
		$$
		Here $C_0$ is the constant given in \eqref{criterion}.
		This is possible because the left side of the first inequality grows like $\delta\sqrt N$ while its right side is bounded by $C\delta$. 
		
		Since for fix $N>0$, $\phi_{1,0}^N$ and $\phi_{2,0}$ satisfy \eqref{potentials} for every $k>3$, then by Theorem~\ref{thm:blow}, the maximal $H^{k+1}$ solution $u^N$ and the corresponding lifespan $T_N^*$ satisfy
		$$
		0<T_N^*
		\leq \frac{8}{\phi_{1,0}^N(0)}
		=\frac{8}{\delta\sqrt{N}\,\rho(0)}
		<\eps.
		$$
		For any $s>k+1$, let $u^{N,s}$ be the maximal $H^s$
		solution with initial datum $u_0^N$, and denote its maximal
		existence time by $T_{N,s}^*$. By uniqueness and the embedding
		$H^s\hookrightarrow H^{k+1}$,
		$$
		T_{N,s}^*\leq T_N^*,
		\qquad
		u^{N,s}(t)=u^N(t)
		\quad\text{for }0\leq t<T_{N,s}^*.
		$$
		Suppose that $T_{N,s}^*<T_N^*$. Then
		$u^N\in C([0,T_{N,s}^*];H^{k+1}(\R^2))$. Write $S^N=(\nabla u^N+(\nabla u^N)^T)/2$, by $H^{k+1}(\R^2)\hookrightarrow W^{1,\infty}(\R^2)$, we obtain
		$$
		\begin{aligned}
			\int_0^{T_{N,s}^*}
			\|S^N(t)\|_{\dot B^0_{\infty,\infty}}\,\dd t
			&\leq C\int_0^{T_{N,s}^*}
			\|\nabla u^N(t)\|_{L^\infty}\,\dd t\\
			&\leq C T_{N,s}^*
			\sup_{0\leq t\leq T_{N,s}^*}\|u^N(t)\|_{H^{k+1}}
			<\infty.
		\end{aligned}
		$$
		This contradicts Lemma~\ref{lem:blow up criteria} at Sobolev
		order $s$. Hence $T_{N,s}^*=T_N^*$ for every $s>4$.
		In particular, for every $0<T<T_N^*$,
		$$
		u^N\in
		\bigcap_{\substack{s>4}}
		C([0,T];H^s(\R^2)).
		$$
		By Bernstein's inequality,
		$$\|S^N\|_{\dot B^0_{\infty,\infty}}\leq C\|u^N\|_{H^2}.$$ 
		If $\sup_{t<T_N^*}\|u^N(t)\|_{H^2}<\infty$, combining with $T_N^*<\infty$, we have
		$$
		\int_0^{T_N^*}\|S^N(t)\|_{\dot B^0_{\infty,\infty}}\,\dd t
		\leq C T_N^*\sup_{t<T_N^*}\|u^N(t)\|_{H^2}<\infty,
		$$
		contradicting Lemma~\ref{lem:blow up criteria}. Therefore, the $H^2$ norm is unbounded on $[0,T_N^*)$, and hence there exists a time $0<t_\eps<T_N^*<\eps$ such that $\|u^N(t_\eps)\|_{H^2}>\eps^{-1}$.
	\end{proof}
	
	\section{Discussion}\label{sec:discussion}
	
	The blow-up criterion in Theorem~\ref{thm:blow} and the ill-posedness result in Theorem~\ref{thm:inflation} are both based on the div--curl decomposition with radial potentials in two dimensions. Building on this structure, we discuss possible directions for further investigation and the technical obstructions that may arise. We first consider global existence conditions for the two-dimensional Euler--Poincar\'e equations within this structure. Once a tangential momentum component is introduced, controlling the Riccati equation at the origin requires additional estimates, and several obstructions arise in deriving global existence conditions. We then discuss further questions that merit investigation and some possible approaches to them.
	
	\paragraph{Global existence in two dimensions.}
	When $q_0=0$, the tangential momentum remains zero. Under the regularity and low-frequency hypotheses of \cite{LYZ2013}, the corresponding radial gradient-momentum solution is global if $\phi_{1,0}(0)\leq 0$ and $p_0(r)=\partial_r\phi_{1,0}(r)\geq 0$. The sign condition on $p_0$ is propagated by the homogeneous transport equation satisfied by the radial momentum. Since $p=\partial_r\phi_1\geq0$ and $\phi_1$ tends to zero at infinity, we have $\phi_1(t,0)\leq\phi_1(t,r)\leq0$. In this case, \eqref{phi_1(t,0)-eq} also yields $\partial_t\phi_1(t,0)\geq0$. Thus the initial central value controls $\|\phi_1(t)\|_{L^\infty}$, and Lemma~\ref{lem:continuation} then guarantees global existence of the solution.
	
	For $q_0\ne0$, let $R(t,a)$ be the radial characteristic defined in the proof of Lemma~\ref{lem:core}, and
	$$
	D(t,a)=\exp\left(\int_0^t
	\left(2v_r+\frac{v}{r}\right)(s,R(s,a))\,\dd s\right).
	$$
	The first equation in \eqref{p q-eq} gives
	$$
	D(t,a)p(t,R(t,a))
	=p_0(a)-\int_0^t D(s,a)
	\left(w_r-\frac{w}{r}\right)(s,R(s,a))q(s,R(s,a))\,\dd s .
	$$
	Consequently, suppose that $p_0(a)\geq0$ for all $a>0$. Since $D(t,a)>0$, to ensure that $p(t,r)\geq0$ throughout the evolution, the second term on the right-hand side of the preceding identity, taken without the preceding minus sign, must not exceed $p_0(a)$ for every $a>0$ and every $0\leq t<T^*$. The second equation in \eqref{p q-eq} preserves the sign of $q$ along each characteristic, but the sign of $(w_r-w/r)q$ also depends on that of $w_r-w/r$.
	
	Moreover, if we can ensure that $p(t,r)\geq 0$ throughout the evolution, then $\|\phi_1(t)\|_{L^\infty}=-\phi_1(t,0)$. To apply the continuation criterion, we therefore need a lower bound on $-\phi_1(t,0)$, i.e., an upper bound for $\phi_1(t,0)$. Under this nonnegativity condition, the identity \eqref{phi_1(t,0)-eq} gives
	$$
	\|\phi_1(t)\|_{L^\infty}
	=-\phi_1(t,0)\leq-\phi_{1,0}(0)+\int_0^t I(s)\,\dd s .
	$$
	Thus, if $\int_0^t I(s)\,\dd s$ remains bounded on every finite time interval and $p\geq0$ is propagated, then Lemma~\ref{lem:continuation} would yield global existence. However, in the proof of Theorem~\ref{thm:blow}, the energy estimate for $I$ relies on the support propagation time of $q$ and therefore yields a bound, determined by the initial data and the conserved energy $E$, only on a finite time interval. Besides, $\|\phi_1(t)\|_{L^\infty}$ cannot be controlled by the $L^2$ growth estimate \eqref{phi1 growth} near the origin. Therefore, a global existence argument for nonzero tangential momentum would require additional and possibly more delicate estimates.
	
	\paragraph{The three-dimensional case.} 
	In three dimensions, we may assume similarly
	$$
	m=\nabla\phi_1+\nabla\times (\phi_2e_z),\quad e_z=(0,0,1),
	$$
	where $\phi_1$ and $\phi_2$ are radial. Naturally, we rewrite the three-dimensional equations in cylindrical coordinates $(\rho,\theta,z)$ with radial potentials, 
	$$
	m=\nabla\phi_1+\nabla\times (\phi_2e_z)=\frac 1 r(\rho \phi'(r) e_\rho-\rho\phi_2'(r)e_\theta+z\phi_1'(r)e_z)
	$$
	with
	$$
	\rho=\sqrt{x^2+y^2},\quad r=\sqrt{x^2+y^2+z^2},\quad e_\rho=(\cos\theta,\sin\theta,0),\quad e_\theta=(-\sin\theta,\cos\theta,0).$$
	Then $m$ and $u$ can be written as
	$$
	m=m^\rho(t,\rho,z)e_\rho+m^\theta(t,\rho,z)e_\theta+m^z(t,\rho,z)e_z,
	$$
	$$u=u^\rho(t,\rho,z)e_\rho+u^\theta(t,\rho,z)e_\theta+u^z(t,\rho,z)e_z.$$
	In particular, $m^\theta e_\theta$ is divergence-free and $\rho m^\theta$ satisfies
	$$
	\left(\partial_t+u^\rho\partial_\rho+u^z\partial_z\right)(\rho m^\theta)
	+(\diver u)\,\rho m^\theta=0.
	$$
	Hence the support of $\rho m^\theta$ is transported along characteristics in the $(\rho,z)$ plane. However, the main difficulty is that one can not show that potentials is  radial  even if the initial potentials are assumed to be radial. The resulting evolution equation becomes more complicated and requires energy estimates for additional quantities, which constitutes another major difficulty.
	
	\paragraph{Nonradial potentials.}
	It is natural to ask whether Theorem~\ref{thm:blow} remains valid for
	$$
	m_0=\nabla\phi_{1,0}+\nabla^\perp\phi_{2,0}
	$$
	when $\phi_{1,0}$ and $\phi_{2,0}$ are not radial in two dimensions. Consider the flow map
	$$
	X'(t)=u(t,X(t)),\qquad X(0)=a,
	$$
	and since $\diver u=\phi_1-(1-\Delta)^{-1}\phi_1$, then along this trajectory we deduce that
	$$
	(\diver u)(t,X(t,a))
	=\bigl[(1-\Delta)^{-1}\phi_1-\phi_1\bigr](t,X(t,a)).
	$$
	A possible strategy is to choose $a$ such that $(\diver u)(t,X(t,a))$ satisfies a Riccati-type inequality that leads to blow-up. For nonradial data, the difficulty is to control the coupling terms between the gradient and divergence-free components along $X(t,a)$. Estimates for these terms would yield conditions on the initial potentials that ensure finite-time blow-up without radial symmetry.
	
	\paragraph{Global weak solutions.}
	In view of the natural $H^1$ energy conservation law \eqref{conservation laws}, it is natural to ask whether the equations admit global weak solutions. More precisely, one seeks a distributional solution $u\in L^\infty(0,\infty;H^1(\mathbb R^2))$, under weaker assumptions on the initial data, satisfying the energy inequality
	$$
	\|u(t)\|_{H^1}^2\leq \|u_0\|_{H^1}^2
	\qquad\text{for almost every }t>0.
	$$
	One possible approach is to start from the global strong solutions in \cite{LYZ2013} and pass to the limit by a compactness argument. The main difficulty is that the $H^1$ bound on the approximation sequence $u_n$ yields only weak compactness of $\nabla u_n$ in $L^2_{\mathrm{loc}}$, which does not determine the weak limit of quadratic products such as $\partial_i u_n\cdot\partial_j u_n$. One would therefore need either additional compactness for the velocity gradients or a structural argument exploiting the equation to show that the defect measure arising in the limit of the quadratic terms vanishes.

	\smallskip
	\noindent\textbf{Acknowledgments}
	This work was partially supported by the National Key R\&D Program of China(No. 2021YFA1002100) and the National Natural Science Foundation of China ( No. 12471233 and No. 12571261).
	
	\smallskip
	\noindent\textbf{Data availability statement}
	Data sharing is not applicable to this article as no datasets were generated or analysed during
	the current study.
	
	\smallskip
	\noindent\textbf{Conflict of interest} On behalf of all authors, the corresponding author states that there are no competing interests. All authors read and approved the final manuscript.
	
	

	\phantomsection
	\addcontentsline{toc}{section}{\refname}
	\bibliographystyle{abbrv} 
	\bibliography{EPrefs}

@article{CL2012,
  author={Chae, Dongho and Liu, Jian-Guo},
  title={Blow-up, zero $\alpha$ limit and the {Liouville} type theorem for the {Euler--Poincar\'e} equations},
  journal={Communications in Mathematical Physics},
  volume={314}, pages={671--687}, year={2012},
  doi={10.1007/s00220-012-1534-8}
}

@article{LYZ2013,
  author={Li, Dong and Yu, Xinwei and Zhai, Zhichun},
  title={On the {Euler--Poincar\'e} equation with non-zero dispersion},
  journal={Archive for Rational Mechanics and Analysis},
  volume={210}, pages={955--974}, year={2013},
  doi={10.1007/s00205-013-0662-4}
}

@inproceedings{HMA2001,
	author    = {Hirani, Anil N. and Marsden, Jerrold E. and Arvo, James},
	title     = {Averaged Template Matching Equations},
	booktitle = {Energy Minimization Methods in Computer Vision and Pattern Recognition},
	series    = {Lecture Notes in Computer Science},
	volume    = {2134},
	pages     = {528--543},
	publisher = {Springer},
	year      = {2001},
	doi       = {10.1007/3-540-44745-8_35}
}

@incollection {HM2005,
	AUTHOR = {Holm, Darryl D. and Marsden, Jerrold E.},
	TITLE = {Momentum maps and measure-valued solutions (peakons,
	filaments, and sheets) for the {EPD}iff equation},
	BOOKTITLE = {The breadth of symplectic and {P}oisson geometry},
	SERIES = {Progr. Math.},
	VOLUME = {232},
	PAGES = {203--235},
	PUBLISHER = {Birkh\"auser Boston, Boston, MA},
	YEAR = {2005},
	ISBN = {0-8176-3565-3},
	MRCLASS = {37K65 (35Q51 37K05 37K40 58D05)},
	MRNUMBER = {2103008},
	MRREVIEWER = {Boris\ A.\ Khesin},
	DOI = {10.1007/0-8176-4419-9\_8},
	URL = {https://doi.org/10.1007/0-8176-4419-9_8},
}

@book {HS2009,
	AUTHOR = {Holm, Darryl D. and Schmah, Tanya and Stoica, Cristina},
	TITLE = {Geometric mechanics and symmetry},
	SERIES = {Oxford Texts in Applied and Engineering Mathematics},
	VOLUME = {12},
	NOTE = {From finite to infinite dimensions,
	With solutions to selected exercises by David C. P. Ellis},
	PUBLISHER = {Oxford University Press, Oxford},
	YEAR = {2009},
	PAGES = {xvi+515},
	ISBN = {978-0-19-921291-0},
	MRCLASS = {37-01 (37Jxx 37K05 37K65 70G45 70G65 70H05 70Hxx)},
	MRNUMBER = {2548736},
	MRREVIEWER = {Frans\ Cantrijn},
}

@article {YY2015,
	AUTHOR = {Yan, Kai and Yin, Zhaoyang},
	TITLE = {On the initial value problem for higher dimensional
	{C}amassa-{H}olm equations},
	JOURNAL = {Discrete Contin. Dyn. Syst.},
	FJOURNAL = {Discrete and Continuous Dynamical Systems},
	VOLUME = {35},
	YEAR = {2015},
	NUMBER = {3},
	PAGES = {1327--1358},
	ISSN = {1078-0947,1553-5231},
	MRCLASS = {35Q53 (35B30 35B44 35D35 35G25)},
	MRNUMBER = {3277196},
	MRREVIEWER = {Amin\ Esfahani},
	DOI = {10.3934/dcds.2015.35.1327},
	URL = {https://doi.org/10.3934/dcds.2015.35.1327},
}

@article {LY2020,
	AUTHOR = {Luo, Wei and Yin, Zhaoyang},
	TITLE = {Blow-up phenomena, ill-posedness and peakon solutions for the
	periodic {E}uler-{P}oincar\'e{} equations},
	JOURNAL = {J. Differential Equations},
	FJOURNAL = {Journal of Differential Equations},
	VOLUME = {268},
	YEAR = {2020},
	NUMBER = {4},
	PAGES = {1307--1325},
	ISSN = {0022-0396,1090-2732},
	MRCLASS = {35Q53 (35B30 35B44 35C07 35G25)},
	MRNUMBER = {4042343},
	DOI = {10.1016/j.jde.2019.08.042},
	URL = {https://doi.org/10.1016/j.jde.2019.08.042},
}

@book {BCD,
	AUTHOR = {Bahouri, Hajer and Chemin, Jean-Yves and Danchin, Rapha\"{e}l},
	TITLE = {Fourier analysis and nonlinear partial differential equations},
	SERIES = {Grundlehren der mathematischen Wissenschaften [Fundamental
	Principles of Mathematical Sciences]},
	VOLUME = {343},
	PUBLISHER = {Springer, Heidelberg},
	YEAR = {2011},
	PAGES = {xvi+523},
	ISBN = {978-3-642-16829-1},
	MRCLASS = {35-02 (35L72 35Q30 42-02 42B37 76B03 76D03 76N10)},
	MRNUMBER = {2768550},
	MRREVIEWER = {Peter\ R.\ Massopust},
	DOI = {10.1007/978-3-642-16830-7},
	URL = {https://doi.org/10.1007/978-3-642-16830-7},
}

@book{Evans2010,
	author    = {Evans, Lawrence C.},
	title     = {Partial Differential Equations},
	edition   = {2},
	series    = {Graduate Studies in Mathematics},
	volume    = {19},
	publisher = {American Mathematical Society},
	address   = {Providence, RI},
	year      = {2010}
}

@incollection {HMR2002,
	AUTHOR = {Holm, Darryl D. and Marsden, Jerrold E. and Ratiu, Tudor S.},
	TITLE = {The {E}uler-{P}oincar\'e{} equations in geophysical fluid
	dynamics},
	BOOKTITLE = {Large-scale atmosphere-ocean dynamics, {V}ol.\ {II}},
	PAGES = {251--300},
	PUBLISHER = {Cambridge Univ. Press, Cambridge},
	YEAR = {2002},
	ISBN = {0-521-80757-3},
	MRCLASS = {37N10 (76B99 76M30 76M60 76U05 86A05)},
	MRNUMBER = {1948165},
	MRREVIEWER = {M.\ E.\ Sousa-Dias},
}

@article {HMR1998,
	AUTHOR = {Holm, Darryl D. and Marsden, Jerrold E. and Ratiu, Tudor S.},
	TITLE = {The {E}uler-{P}oincar\'e{} equations and semidirect products
	with applications to continuum theories},
	JOURNAL = {Adv. Math.},
	FJOURNAL = {Advances in Mathematics},
	VOLUME = {137},
	YEAR = {1998},
	NUMBER = {1},
	PAGES = {1--81},
	ISSN = {0001-8708,1090-2082},
	MRCLASS = {58F05 (58E30 58F40 70H30 73B99 73V25)},
	MRNUMBER = {1627802},
	MRREVIEWER = {Jair\ Koiller},
	DOI = {10.1006/aima.1998.1721},
	URL = {https://doi.org/10.1006/aima.1998.1721},
}

@article {HS2003,
	AUTHOR = {Holm, Darryl D. and Staley, Martin F.},
	TITLE = {Wave structure and nonlinear balances in a family of
	evolutionary {PDE}s},
	JOURNAL = {SIAM J. Appl. Dyn. Syst.},
	FJOURNAL = {SIAM Journal on Applied Dynamical Systems},
	VOLUME = {2},
	YEAR = {2003},
	NUMBER = {3},
	PAGES = {323--380},
	ISSN = {1536-0040},
	MRCLASS = {76D33 (35Q35 37N10)},
	MRNUMBER = {2031278},
	DOI = {10.1137/S1111111102410943},
	URL = {https://doi.org/10.1137/S1111111102410943},
}

@article {CH,
	AUTHOR = {Camassa, Roberto and Holm, Darryl D.},
	TITLE = {An integrable shallow water equation with peaked solitons},
	JOURNAL = {Phys. Rev. Lett.},
	FJOURNAL = {Physical Review Letters},
	VOLUME = {71},
	YEAR = {1993},
	NUMBER = {11},
	PAGES = {1661--1664},
	ISSN = {0031-9007,1079-7114},
	MRCLASS = {35Q51 (58F07 76B15 76B25)},
	MRNUMBER = {1234453},
	DOI = {10.1103/PhysRevLett.71.1661},
	URL = {https://doi.org/10.1103/PhysRevLett.71.1661},
}

@article {CH2,
	AUTHOR = {Constantin, Adrian and Lannes, David},
	TITLE = {The hydrodynamical relevance of the {C}amassa-{H}olm and
	{D}egasperis-{P}rocesi equations},
	JOURNAL = {Arch. Ration. Mech. Anal.},
	FJOURNAL = {Archive for Rational Mechanics and Analysis},
	VOLUME = {192},
	YEAR = {2009},
	NUMBER = {1},
	PAGES = {165--186},
	ISSN = {0003-9527,1432-0673},
	MRCLASS = {35Q53 (35R35 74J30 76A02 76B15)},
	MRNUMBER = {2481064},
	MRREVIEWER = {Bengt\ O.\ Enflo},
	DOI = {10.1007/s00205-008-0128-2},
	URL = {https://doi.org/10.1007/s00205-008-0128-2},
}

@article {CH4,
	AUTHOR = {Constantin, Adrian},
	TITLE = {On the scattering problem for the {C}amassa-{H}olm equation},
	JOURNAL = {R. Soc. Lond. Proc. Ser. A Math. Phys. Eng. Sci.},
	FJOURNAL = {The Royal Society of London. Proceedings. Series A.
	Mathematical, Physical and Engineering Sciences},
	VOLUME = {457},
	YEAR = {2001},
	NUMBER = {2008},
	PAGES = {953--970},
	ISSN = {1364-5021,1471-2946},
	MRCLASS = {37K15 (34L25 35Q53)},
	MRNUMBER = {1875310},
	MRREVIEWER = {Dmitry\ E.\ Pelinovsky},
	DOI = {10.1098/rspa.2000.0701},
	URL = {https://doi.org/10.1098/rspa.2000.0701},
}

@article {CH5,
	AUTHOR = {Constantin, A. and McKean, H. P.},
	TITLE = {A shallow water equation on the circle},
	JOURNAL = {Comm. Pure Appl. Math.},
	FJOURNAL = {Communications on Pure and Applied Mathematics},
	VOLUME = {52},
	YEAR = {1999},
	NUMBER = {8},
	PAGES = {949--982},
	ISSN = {0010-3640,1097-0312},
	MRCLASS = {37K10 (34L40 35Q58 37K15)},
	MRNUMBER = {1686969},
	MRREVIEWER = {Piotr\ G.\ Grinevich},
	DOI = {10.1002/(SICI)1097-0312(199908)52:8<949::AID-CPA3>3.0.CO;2-D},
	URL =
	{https://doi.org/10.1002/(SICI)1097-0312(199908)52:8<949::AID-CPA3>3.0.CO;2-D},
}

@article {CH6,
	AUTHOR = {Constantin, Adrian},
	TITLE = {The {H}amiltonian structure of the {C}amassa-{H}olm equation},
	JOURNAL = {Exposition. Math.},
	FJOURNAL = {Expositiones Mathematicae. International Journal},
	VOLUME = {15},
	YEAR = {1997},
	NUMBER = {1},
	PAGES = {53--85},
	ISSN = {0723-0869},
	MRCLASS = {76B15 (35Q53 58F07)},
	MRNUMBER = {1438436},
	MRREVIEWER = {Mark\ D.\ Groves},
}

@article {CH7,
	AUTHOR = {Fuchssteiner, B. and Fokas, A. S.},
	TITLE = {Symplectic structures, their {B}\"{a}cklund transformations
	and hereditary symmetries},
	JOURNAL = {Phys. D},
	FJOURNAL = {Physica D. Nonlinear Phenomena},
	VOLUME = {4},
	YEAR = {1981/82},
	NUMBER = {1},
	PAGES = {47--66},
	ISSN = {0167-2789,1872-8022},
	MRCLASS = {58F05 (35Q20 58F07 58F37)},
	MRNUMBER = {636470},
	MRREVIEWER = {Toru\ Tsujishita},
	DOI = {10.1016/0167-2789(81)90004-X},
	URL = {https://doi.org/10.1016/0167-2789(81)90004-X},
}

@article {CH8,
	AUTHOR = {Constantin, Adrian and Escher, Joachim},
	TITLE = {Global existence and blow-up for a shallow water equation},
	JOURNAL = {Ann. Scuola Norm. Sup. Pisa Cl. Sci. (4)},
	FJOURNAL = {Annali della Scuola Normale Superiore di Pisa. Classe di
	Scienze. Serie IV},
	VOLUME = {26},
	YEAR = {1998},
	NUMBER = {2},
	PAGES = {303--328},
	ISSN = {0391-173X,2036-2145},
	MRCLASS = {35Q53 (35B40 76B15)},
	MRNUMBER = {1631589},
	MRREVIEWER = {Darryl\ D.\ Holm},
	URL = {http://www.numdam.org/item?id=ASNSP_1998_4_26_2_303_0},
}

@article {CH9,
	AUTHOR = {Constantin, Adrian and Escher, Joachim},
	TITLE = {Well-posedness, global existence, and blowup phenomena for a
	periodic quasi-linear hyperbolic equation},
	JOURNAL = {Comm. Pure Appl. Math.},
	FJOURNAL = {Communications on Pure and Applied Mathematics},
	VOLUME = {51},
	YEAR = {1998},
	NUMBER = {5},
	PAGES = {475--504},
	ISSN = {0010-3640,1097-0312},
	MRCLASS = {35Q53 (35B60 35D05 35L67)},
	MRNUMBER = {1604278},
	MRREVIEWER = {Bruno\ Scheurer},
	DOI = {10.1002/(SICI)1097-0312(199805)51:5<475::AID-CPA2>3.0.CO;2-5},
	URL =
	{https://doi.org/10.1002/(SICI)1097-0312(199805)51:5<475::AID-CPA2>3.0.CO;2-5},
}

@article {CH10,
	AUTHOR = {Danchin, Rapha\"{e}l},
	TITLE = {A few remarks on the {C}amassa-{H}olm equation},
	JOURNAL = {Differential Integral Equations},
	FJOURNAL = {Differential and Integral Equations. An International Journal
	for Theory \& Applications},
	VOLUME = {14},
	YEAR = {2001},
	NUMBER = {8},
	PAGES = {953--988},
	ISSN = {0893-4983},
	MRCLASS = {35Q53 (35B30 35B40 35G25 76B03 76B15)},
	MRNUMBER = {1827098},
	MRREVIEWER = {Kenji\ Nakanishi},
}

@article {CH11,
	AUTHOR = {Rodr\'{\i}guez-Blanco, Guillermo},
	TITLE = {On the {C}auchy problem for the {C}amassa-{H}olm equation},
	JOURNAL = {Nonlinear Anal.},
	FJOURNAL = {Nonlinear Analysis. Theory, Methods \& Applications. An
	International Multidisciplinary Journal},
	VOLUME = {46},
	YEAR = {2001},
	NUMBER = {3},
	PAGES = {309--327},
	ISSN = {0362-546X,1873-5215},
	MRCLASS = {35Q53 (35B30 35B40 35K90)},
	MRNUMBER = {1851854},
	MRREVIEWER = {Muthusamy\ Vanninathan},
	DOI = {10.1016/S0362-546X(01)00791-X},
	URL = {https://doi.org/10.1016/S0362-546X(01)00791-X},
}

@article {CH12,
	AUTHOR = {Constantin, Adrian and Escher, Joachim},
	TITLE = {Wave breaking for nonlinear nonlocal shallow water equations},
	JOURNAL = {Acta Math.},
	FJOURNAL = {Acta Mathematica},
	VOLUME = {181},
	YEAR = {1998},
	NUMBER = {2},
	PAGES = {229--243},
	ISSN = {0001-5962,1871-2509},
	MRCLASS = {35Q35 (76B15)},
	MRNUMBER = {1668586},
	MRREVIEWER = {Alan\ Jeffrey},
	DOI = {10.1007/BF02392586},
	URL = {https://doi.org/10.1007/BF02392586},
}

@article {CH13,
	AUTHOR = {Constantin, Adrian},
	TITLE = {On the {C}auchy problem for the periodic {C}amassa-{H}olm
	equation},
	JOURNAL = {J. Differential Equations},
	FJOURNAL = {Journal of Differential Equations},
	VOLUME = {141},
	YEAR = {1997},
	NUMBER = {2},
	PAGES = {218--235},
	ISSN = {0022-0396,1090-2732},
	MRCLASS = {35Q53 (34G20 76B15)},
	MRNUMBER = {1488351},
	MRREVIEWER = {M\'arcia\ Scialom},
	DOI = {10.1006/jdeq.1997.3333},
	URL = {https://doi.org/10.1006/jdeq.1997.3333},
}

@article {CH15,
	AUTHOR = {McKean, H. P.},
	TITLE = {Breakdown of a shallow water equation},
	NOTE = {Mikio Sato: a great Japanese mathematician of the twentieth
	century},
	JOURNAL = {Asian J. Math.},
	FJOURNAL = {Asian Journal of Mathematics},
	VOLUME = {2},
	YEAR = {1998},
	NUMBER = {4},
	PAGES = {867--874},
	ISSN = {1093-6106,1945-0036},
	MRCLASS = {35Q58 (35Q51 37K40)},
	MRNUMBER = {1734131},
	MRREVIEWER = {Alan\ Shapiro},
	DOI = {10.4310/AJM.1998.v2.n4.a10},
	URL = {https://doi.org/10.4310/AJM.1998.v2.n4.a10},
}

@article {CH17,
	AUTHOR = {Constantin, Adrian and Strauss, Walter A.},
	TITLE = {Stability of peakons},
	JOURNAL = {Comm. Pure Appl. Math.},
	FJOURNAL = {Communications on Pure and Applied Mathematics},
	VOLUME = {53},
	YEAR = {2000},
	NUMBER = {5},
	PAGES = {603--610},
	ISSN = {0010-3640,1097-0312},
	MRCLASS = {35Q51 (35B35 35Q53 76B15)},
	MRNUMBER = {1737505},
	MRREVIEWER = {Marcel\ Oliver},
	DOI = {10.1002/(sici)1097-0312(200005)53:5<603::aid-cpa3>3.0.co;2-l},
	URL =
	{https://doi.org/10.1002/(sici)1097-0312(200005)53:5<603::aid-cpa3>3.0.co;2-l},
}

@article{CH18,
	author  = {Danchin, Rapha{\"e}l},
	title   = {A note on well-posedness for {Camassa--Holm} equation},
	journal = {Journal of Differential Equations},
	volume  = {192},
	number  = {2},
	pages   = {429--444},
	year    = {2003},
	doi     = {10.1016/S0022-0396(03)00096-2}
}

@article {CH19,
	AUTHOR = {Constantin, Adrian and Escher, Joachim},
	TITLE = {Global weak solutions for a shallow water equation},
	JOURNAL = {Indiana Univ. Math. J.},
	FJOURNAL = {Indiana University Mathematics Journal},
	VOLUME = {47},
	YEAR = {1998},
	NUMBER = {4},
	PAGES = {1527--1545},
	ISSN = {0022-2518,1943-5258},
	MRCLASS = {35Q53 (35D05 76B15)},
	MRNUMBER = {1687106},
	MRREVIEWER = {Marcel\ Oliver},
	DOI = {10.1512/iumj.1998.47.1466},
	URL = {https://doi.org/10.1512/iumj.1998.47.1466},
}

@article{CH20,
	author  = {Constantin, Adrian and Molinet, Luc},
	title   = {Global weak solutions for a shallow water equation},
	journal = {Communications in Mathematical Physics},
	volume  = {211},
	number  = {1},
	pages   = {45--61},
	year    = {2000},
	doi     = {10.1007/s002200050801}
}

@article {CH21,
	AUTHOR = {Xin, Zhouping and Zhang, Ping},
	TITLE = {On the weak solutions to a shallow water equation},
	JOURNAL = {Comm. Pure Appl. Math.},
	FJOURNAL = {Communications on Pure and Applied Mathematics},
	VOLUME = {53},
	YEAR = {2000},
	NUMBER = {11},
	PAGES = {1411--1433},
	ISSN = {0010-3640,1097-0312},
	MRCLASS = {35Q53 (35D05 76B15)},
	MRNUMBER = {1773414},
	MRREVIEWER = {Sergey\ Nikolaevich\ Alekseenko},
	DOI = {10.1002/1097-0312(200011)53:11<1411::aid-cpa4>3.0.co;2-5},
	URL =
	{https://doi.org/10.1002/1097-0312(200011)53:11<1411::aid-cpa4>3.0.co;2-5},
}

@article {CH22,
	AUTHOR = {Xin, Zhouping and Zhang, Ping},
	TITLE = {On the uniqueness and large time behavior of the weak
	solutions to a shallow water equation},
	JOURNAL = {Comm. Partial Differential Equations},
	FJOURNAL = {Communications in Partial Differential Equations},
	VOLUME = {27},
	YEAR = {2002},
	NUMBER = {9-10},
	PAGES = {1815--1844},
	ISSN = {0360-5302,1532-4133},
	MRCLASS = {35Q53 (35B40 76B03 76B15)},
	MRNUMBER = {1941659},
	DOI = {10.1081/PDE-120016129},
	URL = {https://doi.org/10.1081/PDE-120016129},
}
\end{document}